\documentclass[msom, nonblindrev]{informs3} % current default for manuscript submission

\OneAndAHalfSpacedXI % current default line spacing
\usepackage{preamble}

\TheoremsNumberedThrough     % Preferred (Theorem 1, Lemma 1, Theorem 2)
\EquationsNumberedThrough    % Default: (1), (2), ...
\MANUSCRIPTNO{MSOM-2024-1346}

\begin{document}
%%%%%%%%%%%%%%%%

%\setcounter{page}{1}
%\DOI{10.1287/xxxx.0000.0000}%

\RUNAUTHOR{Wouda, Aerts-Veenstra, and Van Foreest}

\RUNTITLE{An integrated selection and routing policy for urban waste collection}
\TITLE{An integrated selection and routing policy for urban waste collection}


\ARTICLEAUTHORS{%
\AUTHOR{Niels A. Wouda, Marjolein Aerts-Veenstra, Nicky van Foreest}
\AFF{Department of Operations, University of Groningen, \EMAIL{\{n.a.wouda, m.aerts-veenstra, n.d.van.foreest\}@rug.nl}}
} % end of the block

\ABSTRACT{%
\textbf{Problem definition}: 
We study a daily urban waste collection problem arising in the municipality of Groningen, The Netherlands, where residents bring their waste to local underground waste containers organised in clusters.
The municipality plans routes for waste collection vehicles to empty the container clusters.
These routes should be as short as possible to limit operational costs, but also long enough to visit sufficiently many clusters and ensure that containers do not overflow.
A complicating factor is that the actual fill levels of the clusters' containers are not known, and only the number of deposits is observed.
Additionally, it is unclear whether the containers should be upgraded with expensive fill level sensors so that the service level can be improved or routing costs can be reduced.
\textbf{Methodology/results}:
We propose an eﬀicient integrated selection and routing (ISR) policy that jointly optimises the daily cluster selection and routing decisions.
The integration is achieved by first estimating prizes that express the urgency of selecting a cluster to empty, and then solving a prize-collecting vehicle routing problem with time windows and driver breaks to collect these prizes while minimising routing costs.
We use a metaheuristic to solve the prize-collecting vehicle routing problem inside a realistic simulation environment that models the waste collection problem faced by the municipality.
\textbf{Implications}:
We show that solving the daily waste collection problem in this way is very effective, and can lead to substantial cost savings for the municipality in practice, with no reduction in service level.
In particular, by integrating the container selection and routing problems using our ISR policy, routing costs can be reduced by more than 40\% and the fleet size by 25\%.
Finally, with our procedure to estimate the fill level based on the number of deposits, we show that more advanced measuring techniques do not significantly reduce routing costs, and the service level not at all.
These results should generalise to any municipality operating waste collection from central waste containers.
}%

%\KEYWORDS{deterministic inventory theory; infinite linear programming duality; 
%  existence of optimal policies; semi-Markov decision process; cyclic schedule}

% Authors should select keywords to describe their paper’s theoretical and methodological orientation. A list of keywords is available within ScholarOne Manuscripts. Keywords should appear beneath the abstract in the manuscript file.
% \HISTORY{}

\maketitle
% \tableofcontents
%%%%%%%%%%%%%%%%%%%%%%%%%%%%%%%%%%%%%%%%%%%%%%%%%%%%%%%%%%%%%%%%%%%%%%

% Samples of sectioning (and labeling) in TRSC
% NOTE: (1) \section and \subsection do NOT end with a period
%       (2) \subsubsection and lower need end punctuation
%       (3) capitalization is as shown (title style).

\section{Introduction}
\label{sec:introduction}

%- High level probleem; wat de gemeente doet (namelijk selectieprobleem en VRP met geselecteerde containers), en hoe wij dit geintegreerd aanpakken op basis van data.\\
%- Literatuur (afvalprobleem, vrp etc., evt heel kort iets van methodologie maar grootste deel in sectie methodology)\\
%- Contribution 
%\ma{TODO: either in intro or problem formulation explain the PCVRPTW in words, and solve for only one day.}
% \todo{nvf: als begin vind ik dit een prima intro. Er moet nog wel geschaafd worden. Dat lijkt veel werk, maar dat valt wel mee. Je moet vooral de boel in de goede volgorde zetten, zodat er een begrijpelijk verhaal uitkomkt. Als leidend principe: als je iets opschrijft, vraag je af welke vraag je beantwoordt. Als je dat weet, vraag je af of die vraag op die plek wel nuttig is, maw, of een lezer die vraag op die plek in de tekst zal stellen?  Het valt me dat je wel de goede vragen stelt, en de goede antwoorden geeft, alleen de volgorde is niet zo goed. Maar... het is een zeer bemoedigende ontwikkeling: er is een intro!}

Efficient urban waste collection is a fundamental challenge in smart city operations~\citep{hasija_smart_2020, mak_enabling_2022}.
Waste collection vehicles contribute to air pollution emissions in cities~\citep{maimoun_emissions_2013}.
Additionally, these vehicles increase the traffic burden on crowded city streets~\citep{han_optimizing_2024}, and can be dangerous to other road users like pedestrians and cyclists, particularly during busy hours.
As nearly 80\% of the global population is expected to live in urban environments by 2050~\citep{oecd_cities_2020}, the design and operation of efficient urban waste collection systems will become even more relevant in the near future than they already are today.

\begin{figure}
\centering
\begin{minipage}{.49\textwidth}
    \centering
    \includegraphics[width=\linewidth]{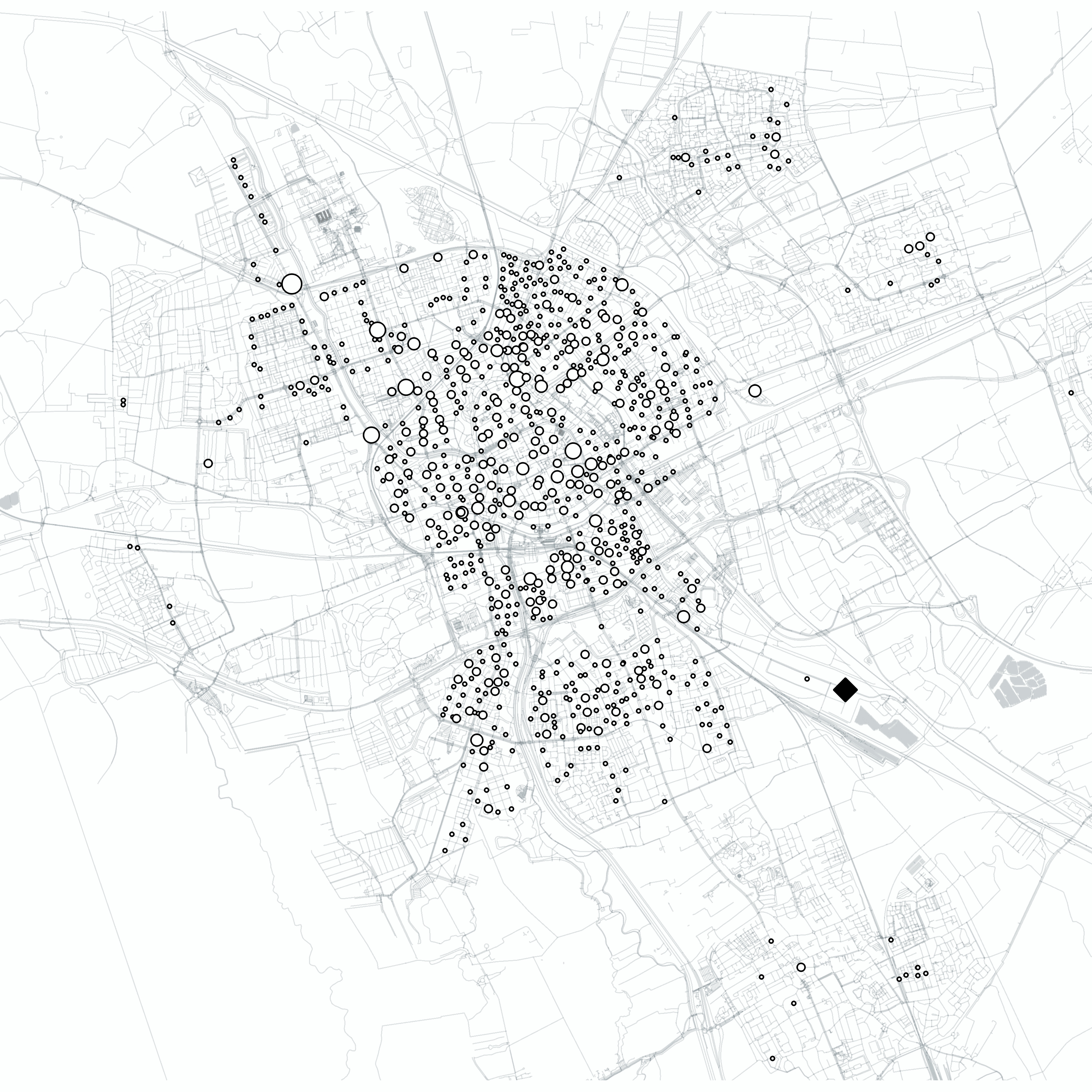}
    \caption{%
        An overview of the clusters (white circles, their size indicates the number of containers) and the depot (black diamond).
    }
    \label{fig:area_of_operations}
\end{minipage}%
\hfill
\begin{minipage}{.49\textwidth}
    \centering
    \includegraphics[width=\linewidth]{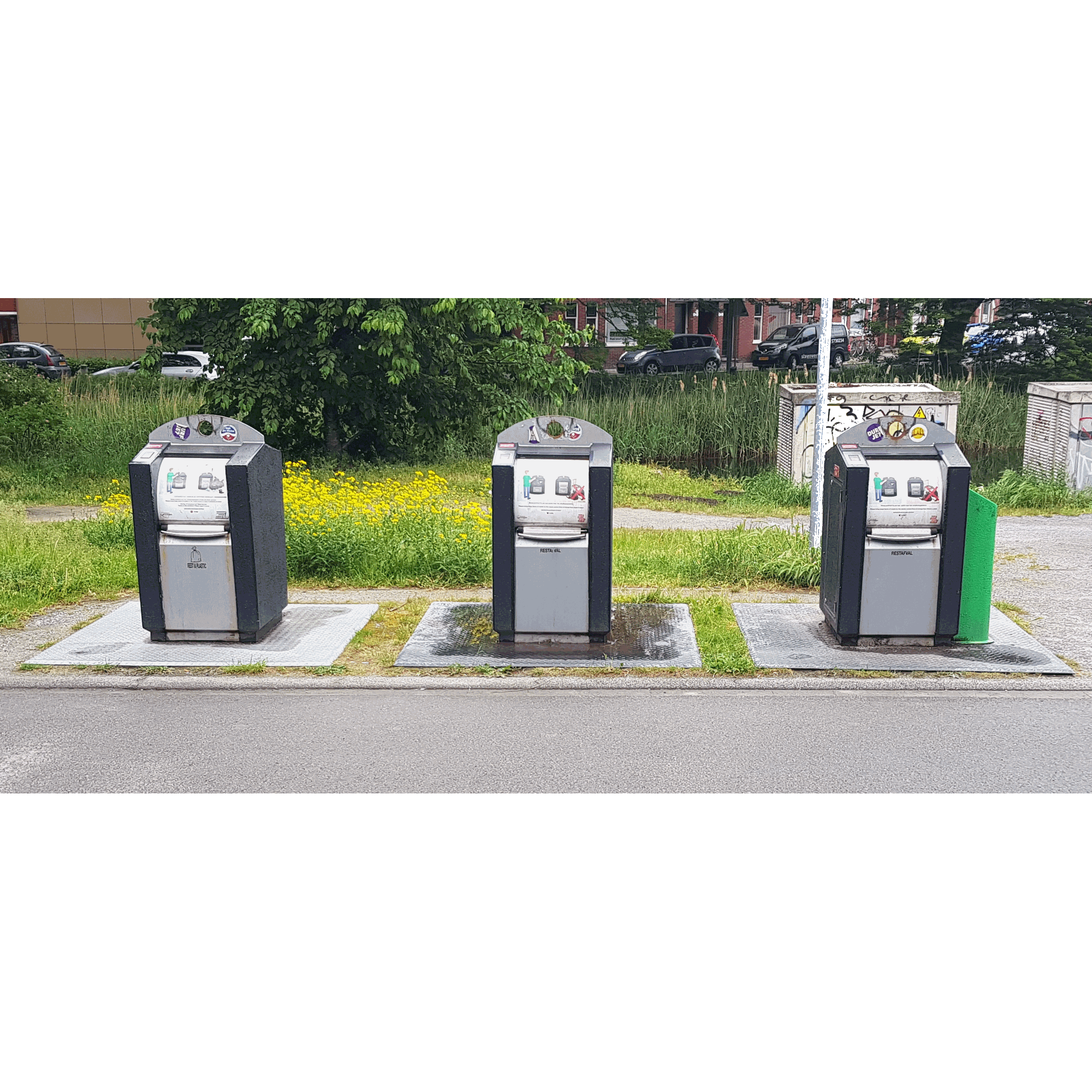}
    \caption{%
        A cluster of three containers.
        Below the sheet metal flooring, each container's underground bin extends down several metres.
    }
    \label{fig:container}
\end{minipage}
\end{figure}

In this paper, we study the operation of an urban waste collection system arising in the municipality of Groningen, The Netherlands, a city with about 240,000 residents on an area of 186 km$^2$.
The majority of these residents use 850 underground container clusters (pictured in~\cref{fig:area_of_operations,fig:container}) placed throughout the city to deposit their household waste.
Each household is assigned to one of these container clusters, typically the one closest to their place of residence.
A container cluster consists of one or more actual containers.
Each household is provided with a card they can use to access any container in their assigned cluster.
The use of these cards is registered in real-time, but the volume and weight of the household waste deposited into the container are not measured.

The municipality is responsible for emptying the clusters in time to avoid overflows, which occur when the cluster is full and newly arriving deposits are placed next to the full containers.
The municipality plans which clusters to empty during daily shifts by planning routes visiting a subset of the container clusters, using the real-time data observed from the cards used to access containers.
In line with the routing literature, we will refer to the emptying of a container cluster as a service.
Since the arrival rates of new deposits vary substantially between clusters, a policy based on fixed routes may service clusters that contain few deposits, while other clusters may have been full for several days before a service is performed~\citep{mcleod_improving_2014}.
Therefore, the plans are made as late as possible, just minutes before the actual shift commences, to ensure the latest card data is used to inform the cluster selection and routing decisions.

The municipality's objective is to minimise overflows over a longer time horizon, while also keeping operational costs related to the routing decisions under control.
This is a challenging problem because the decision of which container clusters to empty today interact with the decisions made tomorrow and beyond: most clusters that are emptied today will not need to be emptied tomorrow, however, when a half-full cluster lies near to a cluster that is about to overflow today, it might still be interesting to empty the half-full cluster today as well, because that consolidation leads to shorter route lengths on later days.
Thus, an objective function for the daily shift planning problem should account for these long-term costs.

We propose an integrated selection and routing (ISR) policy to solve the daily shift planning problem while taking these long-term goals into account.
Our policy constructs a shift plan by first determining a prize value for each container cluster that captures the likelihood the cluster will overflow before the next planning moment, and then solving a prize-collecting vehicle routing problem parametrised by these prizes.
The goal of the routing problem is to minimise driving distance while maximising the value of prizes collected by servicing container clusters.
Thus, by balancing the prizes well, the municipality can balance operational (routing) costs and overflow probabilities.
The prizes are based on the real-time card data and carefully designed to account for the cost of serving a container cluster today, versus postponing the decision until tomorrow.
Finally, the resulting prize-collecting vehicle routing problem is solved using a metaheuristic to obtain the actual shift plan.

Through a large simulation study based on real-world data from the municipality of Groningen, we show that our ISR policy finds very good routes in less than a minute of computational time.
Compared to the municipality's current policy, our policy maintains the excellent service levels the municipality already achieves, but reduces route lengths by 40\% on average.
Additionally, we show that it may be possible to reduce the fleet of collection vehicles by 25\% from four to three vehicles.
Finally, we observe that measuring the actual fill levels of the containers is not necessary for good performance; exploiting the current information, that is, counting the number of deposits, suffices.
Our policy and results are easy to deploy in the current practice of the municipality, and can be applied in any other municipality that uses containers to collect household waste.

\subsection{Related work}

Two systems are in common use for the collection of residential waste.
The first is a door-to-door system, where households place their waste bins by the side of the street.
This system typically requires a fixed collection schedule to ensure households know when to place their waste bins by the side of the street.
The second system is a collection-point system, where households deposit their waste in central neighbourhood waste containers.
When equipped with sensors registering either the arrival of new deposits or the container's current fill level, a container can be serviced just-in-time by a waste collection vehicle.
A significant advantage of the collection-point system in urban environments is that it frees up valuable space on the side of the street~\citep{glaser_waste_2022}.
Additionally, it enables consolidation opportunities that cannot be achieved through a door-to-door system, thus reducing costs~\citep{bautista_solving_2008, rossit_exact_2020}.
In Europe, for example, although both systems are used extensively~\citep{eea_country_2016}, the collection-point system is commonly found in dense urban environments.

The system we are concerned with in this paper is an existing urban collection-point system, where the collection points consist of clusters of underground containers.
These containers are equipped with sensors that register the arrival of individual waste deposits.
The availability of such sensor data extends the operational settings of~\cite{glaser_introduction_2021} and~\cite{lavigne_memetic_2023} by allowing dynamic servicing based on the observed deposit arrivals.
Dynamic servicing can result in substantial cost savings~\citep{johansson_effect_2006,mcleod_improving_2014} since the operator can prioritise visiting only container clusters that are (likely to be) nearly full, rather than having to visit every cluster according to a fixed schedule.
While studies that consider container-based sensors registering the current fill levels abound~(see, for example~\cite{faccio_waste_2011},~\cite{mes_inventory_2014}, and~\cite{mamun_theoretical_2016}), in our case only the arrival of new deposits is observed, not the actual container fill levels.
This setting is similar to that of~\cite{fadda_waste_2018}, who investigate a system that also tracks deposit arrivals, but not fill levels: they gather data on historical fill levels measured upon the moment of service, and use an average fill level for a given number of deposits to decide which containers to empty during each shift.
Working with such averages, however, ignores that the fill levels depends on the unobserved volume of each deposit.
This results in a system where containers will overflow more often than necessary, or increased routing costs.

Solution approaches for waste collection problems vary.
\cite{belien_municipal_2014} survey the literature on municipal waste collection.
They find that heuristic approaches are amongst the most commonly used to generate fixed route plans, while little emphasis is placed on dynamic problems with repeated decision-making.
Solution approaches for dynamic waste collection problems are discussed in the survey of~\cite{ghiani_operations_2014}, but with a focus on strategic decisions, that is, the expansion and operation of an existing waste collection system over longer time frames.

Our urban waste collection problem is a dynamic vehicle routing problem, so solution methods used to solve those problems should also be very suitable in our problem.
A common technique to solve dynamic vehicle routing problems involves the use of roll-out algorithms, which during each decision epoch make a locally optimal---or, at least, good enough---decision.
Such techniques are very popular in the literature: see, for example, the works of~\cite{psaraftis_dynamic_2016,bertazzi_faster_2018,zhang_solving_2023,lan_iterative_2024} for applications and recent developments.
The epoch decisions are often made using look-ahead policies that sample future trajectories, solve the sampled trajectory, and average the sampling solutions using, for example, consensus functions~\citep{bent_scenario-based_2004} to arrive at a decision to implement now.
A downside to such look-ahead policies is their computational expense: to obtain good solutions, they need to sample many trajectories, and solve a routing problem for each trajectory.
To reduce the computational expense, we will instead directly parameterise the epoch decision problem by appropriately accounting for the future state of the system.
This effectively side-steps the sampling costs inherent in many look-ahead approaches.
The idea goes back at least to the solution approach of~\cite{dror_inventoryrouting_1987}, who use penalty and incentive terms to weigh the future effects of local epoch decisions in an inventory routing problem, and has most recently been used by~\cite{baty_combinatorial_2024}, who use machine learning techniques to assign incentives for their routing problem to dispatch or postpone a routing decision in the current epoch of a dynamic same-day delivery problem.

\subsection{Outline}
The rest of this paper is structured as follows.
In~\cref{sec:model_setting} we describe the daily shift planning problem for the waste collection problem in detail.
\Cref{sec:policies} presents our new ISR policy, as well as a baseline policy currently used by the municipality.
\Cref{sec:case_study} forms the core of the paper, where we investigate the performance of our ISR policy through a simulation case study of the urban waste collection system in the municipality of Groningen, and consider several scenarios that are relevant to the municipality.
In~\cref{sec:discussion} we discuss the implications of our results for practice, and reflect on some limitations of our case study approach.
\Cref{sec:conclusion} concludes the paper.

\section{Model setting and framework}
\label{sec:model_setting}

In this section we describe the model setting of our daily waste collection problem.
While we particularly focus on the context of the municipality of Groningen, our description covers a wide range of waste collection systems.

\paragraph{The system}
The municipality of Groningen is responsible for the operation of 850 container clusters, consisting of 1120 individual underground containers.
An overview of the container locations is presented in~\cref{fig:area_of_operations}.
These underground containers have an above-ground drum, shown in~\cref{fig:container}, that allows waste deposits with a maximal volume of 100 litres.
The underground bin has a capacity between 4000 and 6000 litres, depending on the container type.
Each deposit arrival is registered in real-time, but the actual fill levels are not observed.
As deposits accumulate over time, the containers in each cluster fill up, and need to be emptied in time by a collection vehicle \textit{servicing} the container cluster.
When a cluster is visited, the vehicle always services all containers of the cluster.
When timely service does not occur, newly arriving deposits are often placed next to the full container---this is not allowed, but occurs as people want to get rid of their waste.
Clearly, such \textit{overflows} cause significant nuisance.
We note that almost all clusters fill up within a week if unserviced, which means that the duration waste spends inside the cluster is not a factor of concern. 
The service duration at each container cluster consists of two minutes for set-up and tear-down time for the collection vehicle, and one additional minute of handling time for each individual container in the cluster.
For example, a cluster of five containers requires seven minutes of service time.
We remark in passing that drivers can deal with small amounts of overflow, for example, two or three garbage bags, but when the overflow volume is larger, a specialised vehicle has to clean up the excess waste.
Thus, good routing plans should take into consideration the volume of overflow besides just its occurrence.
Finally, although most of the clusters may be serviced at any moment during the day, local legislation requires that the 57 clusters in the inner city are serviced in the morning, before noon.
This requires us to take into account time windows.

The municipality operates from a single depot in the south of Groningen, shown in~\cref{fig:area_of_operations}, where the daily shifts start at 7:00 in the morning.
Four collection vehicles are available for residential waste collection.
Drivers of these vehicles take two breaks at the depot during their shifts: first a coffee break of thirty minutes at 10:00 in the morning, and then a thirty minute lunch break at 12:00.
Neither break is required if the shift concludes before the break is scheduled.
A complete daily shift should ideally not take more than seven hours in total, including the two breaks.
During the breaks, the drivers empty their vehicles at a processing location right next to the depot.
Since the vehicles are regularly emptied, vehicle capacities are not exceeded.

\paragraph{Demand}
Let $C$ denote the set of all container clusters.
The arrival of deposits at each container cluster $c \in C$ is governed by a non-homogeneous Poisson process $N_c(t)$ with known rate function $\lambda_c(t)$.
Historical arrival data provides an excellent means to estimate this rate function.
We assume that the volumes of deposits at cluster $c$ form a sequence of independent random variables $\{ v_{c~i} \}_i$ with common mean $\mu_c > 0$ and variance $\sigma_c^2 > 0$.
These volumes are unobserved, and thus the actual $\mu_c$ and $\sigma_c^2$ are unknown.

\paragraph{Shift planning}
The planner plans a shift on each day, at fixed times $T_1, T_2, \ldots$.
When planning shift $i$ at time $T_i$, the planner prioritises clusters based on the cluster capacity and the number of deposits that arrived since the last moment of service.
We model the priority as a vector of prizes $p \in \mathbb{R}_{\ge 0}^{|C|}$.
In~\cref{sec:policies} we discuss some policies to determine these prizes.

Once the prize vector $p$ is determined, the planner creates the shift plan at time $T_i$ by solving a prize-collecting vehicle routing problem denoted as $\text{VRP}(p)$.
Our formulation of this problem is based on the widely used three-index formulation for the vehicle routing problem with time windows~\citep{golden_implementing_1977, toth_vehicle_2002}.
Let $G$ be a directed graph $G = (N, A)$, where $N$ is the set of nodes and $A \subset N \times N$ is the set of arcs.
We duplicate the depot $0$ as an origin depot $0^o$, and a destination depot $0^d$, each equipped with the same data as the original depot $0$.
We include the set $B$ of driver breaks to be had at the depot location $0$.
Then, $N = C \cup B \cup \{ 0^o, 0^d \}$.
Each node $i \in N$ has a service (respectively, break) duration $s_i \ge 0$, and a time window $[e_i, l_i]$, where $e_i (l_i)$ is the earliest (latest) time at which a service or break can start.
Arriving before the start of a time window incurs waiting time, and arriving after the end of a time window is not allowed.
Each arc $(i, j) \in A$ is associated with a travel cost $d_{ij} \ge 0$ and a travel duration $\tau_{ij} \ge 0$.

Let $x_{ij}^k$ be a binary variable equal to one if arc $(i, j) \in A$ is traversed by vehicle $k \in K$, and 0 otherwise.
Let $y_c^k$ be a binary variable equal to one if cluster $c \in C$ is visited by vehicle $k \in K$, and 0 otherwise.
The continuous variable $\alpha_i^k \ge 0$ determines the time at which service starts at node $i \in N$ by vehicle $k \in K$. 
Our formulation is as follows:
\begin{subequations}
    \label{eq:exact_vrp_formulation}
    \begin{align}
        \text{VRP}(p): \qquad
        \min_{x, y, \alpha} \quad & \sum_{k\in K} \sum_{(i, j) \in A} d_{ij} x_{ij}^k + \sum_{k\in K}\sum_{c \in C} p_c (1 - y_c^k) \label{eq:obj}\\
        \text{s. t.} \quad 
        & \sum_{k \in K} y_c^k \le 1 && \forall c \in C \label{eq:optional} \\
        & \sum_{j: (c, j) \in A} x_{cj}^k = y_c^k  && \forall c \in C, \forall k \in K \label{eq:visits} \\
        & \sum_{b: (b, j) \in A} x_{bj}^k = 1 && \forall b \in B, \forall k \in K \label{eq:breaks} \\
        & \sum_{j: (i,j) \in A} x_{ij}^k - \sum_{j: (j, i) \in A} x_{ji}^k = 
            \begin{cases}
            0 & \text{ if } i \in C \\
            1 & \text{ if } i = 0^o  \\
            -1 & \text{ if } i = 0^d \\
            \end{cases}
        && \forall i \in N, \forall k \in K \label{eq:flow_balance} \\
        & x_{ij}^k = 1 \implies \alpha_{i}^k + s_i + \tau_{ij} \leq \alpha_j^k && \forall (i,j)\in A, \forall k \in K \label{eq:arrival_times}\\
        & \alpha_i^k \in [e_i, l_i] && \forall i \in N, \forall k\in K \label{eq:dom1} \\
        & x_{ij}^k \in \{ 0, 1 \} && \forall (i, j) \in A, \forall k \in K \label{eq:dom2} \\
        & y_c^k \in \{ 0, 1 \} && \forall c \in C, \forall k \in K \label{eq:dom3}
    \end{align}
\end{subequations}
The objective \eqref{eq:obj} is to minimise the travel distance plus the sum of uncollected prizes.
We will sometimes use $p_c = \infty$ to model a required cluster visit.
If, after solving, this results in an optimal objective value of $\infty$, the model should be interpreted as infeasible.
Constraints \eqref{eq:optional} ensure that each container cluster is visited at most once.
Constraints \eqref{eq:visits} couple the binary decision variables $x$ and $y$.
Constraints \eqref{eq:breaks} ensure the drivers of each vehicle $k \in K$ take their required breaks.
Constraints \eqref{eq:flow_balance} are the flow conservation constraints, which ensure routes start and end at the depot.
The implication constraints~\eqref{eq:arrival_times} ensure that if an arc $(i, j) \in A$ is traversed, the arrival time at node $j$ is at least the arrival time at the previously visited node $i$ plus the service time of node $i$ and the travel time of arc $(i,j)$.
Together with constraints~\eqref{eq:dom1}, these ensure that the time windows are respected.
Finally, constraints~\eqref{eq:dom2} and~\eqref{eq:dom3} define the domains of $x$ and $y$, respectively.

Let $S_c(t)$ denote the service counting process of cluster $c \in C$.
Solving the $\text{VRP}(p)$ of~\eqref{eq:exact_vrp_formulation} yields the services planned for shift $i$ as follows.
Let $x^*, y^*, \alpha^*$ be a solution obtained from solving $\text{VRP}(p)$.
For each cluster $c \in C$, when $y^*_c = 1$, we increment the counting process $S_c$ by one at time $\alpha^*_c$.
When $y^*_c = 0$, the cluster is not serviced during this shift, and $S_c$ remains unchanged.

\section{Policies}
\label{sec:policies}

In this section we present two policies to determine the prize vector $p$ to use for shift $i$.
% \todo{Not a big fan of $A_c$ for time of last service, but let's hash out notation for service/emptying next}
Throughout this section, we write
\[ l_c(a, b) = \int_a^b \lambda_c(t) \,\mathrm{d}t,\]
for the expected number of deposits at $c$ in the interval $[a, b$].
Let $A_c(t)$ be the last time that cluster $c \in C$ has been serviced before time $t$.
With this, we write $n_c = N_c(T_i) - N_c(A_c(T_i))$ for the number of deposits that arrived between the last visit of the cluster and the planning moment $T_{i}$.

\Cref{subsec:baseline_policy} describes a baseline policy that the municipality currently uses.
Then, in~\cref{subsec:prize_collecting_policy} we present our ISR policy.
Both policies uses version 0.7.0 of the PyVRP software package~\citep{wouda_pyvrp_2024} to solve the resulting $\text{VRP}(p)$.
Appendix~\ref{app:pyvrp_details} provides details of our modifications to this package to efficiently solve the prize-collecting vehicle routing problem.

\subsection{A baseline policy}
\label{subsec:baseline_policy}

An intuitive policy is to decompose shift planning into two parts: first, select a subset of containers to visit, and then solve the resulting routing problem with just those selected containers.
In fact, the municipality selects 250 container clusters that are expected to fill up the soonest, and then determines a route plan to visit these 250 containers.
Let us explain this \textit{baseline} policy.

The municipality assumes all deposits have the same volume $D = 60~\text{litres}$; they absorb the variability in volume in a correction factor $r_{c}$ for cluster $c \in C$, see below.
Under this assumption, a cluster $c$ can receive $V_c / D$ deposits before it is full.
At time $T_{i}$, the number of deposits is $n_c$, hence a volume of $V_c / D - n_c$ remains for additional deposits after $T_{i}$.
Thus, the expected moment the cluster $c$ is full can be found as
\begin{equation}
    \label{eq:time_till_full_baseline}
    %h_c = \max \left\{ t~\middle|~t l_c(0, T_i) \le r_c \left( \frac{V_c}{D} - n_c \right) \right\},
    h_c = \max \left\{ t~\middle|~ l_c(T_i, t) \le r_c ( V_c/D - n_c ) \right\},
\end{equation}
where $r_c > 0$ is a correction factor used by the municipality that  ranges between $0.8$ and $5$, depending on the cluster.

Let $\pi(1), \ldots, \pi(|C|)$ be a permutation of the clusters $C$ such that $h_{\pi(1)} \le h_{\pi(2)} \le \cdots \le h_{\pi(|C|)}$.  
Then, we set the prize for cluster $c$ as
\[
    p_c =
    \begin{cases}
        \infty & \text{if}~h_c \le h_{\pi(250)}, \\
        0 & \text{otherwise}.
    \end{cases}
\]
These prizes effectively model two sets of clusters: those whose service is required, with prizes set to $\infty$, and those whose service can be skipped, with prizes set to $0$.

Although in the baseline policy drivers take two breaks during their shift, the (black box) routing solve that the municipality currently uses does not take breaks into account explicitly.
Instead, the drivers choose appropriate break moments, but this certainly increases routing costs.
To ensure our results are consistent with current practice, we ignore \eqref{eq:breaks} when solving $\text{VRP}(p)$ under the baseline policy.
Instead, the simulation environment that we will explain in~\cref{subsec:simulator} plans the breaks.

\subsection{An integrated selection and routing policy}
\label{subsec:prize_collecting_policy}

We now present our integrated selection and routing (ISR) policy.
Below we characterise a random variable $Z_{c~i}$ to model the fill level of cluster $c \in C$ at $T_{i + 1}$ when planning the shift at $T_i$.
This allows us to use the overflow probability $\text{Pr}( Z_{c~i} > V_c)$ to set the prize as
\begin{equation}
    \label{eq:isr_prizes}
    p_c =
    \begin{cases}
        \infty & \text{if}~\text{Pr}( Z_{c~i} > V_c) \ge 1 - \epsilon, \\
        \rho \text{Pr}(Z_{c~i} > V_c) & \text{otherwise}.
    \end{cases}
\end{equation}
Here, $\epsilon \in [0, 1]$ is a control parameter that ensures a container is always serviced when its estimated overflow probability is sufficiently large, and $\rho > 0$ is a weight parameter that balances the trade-off between distance cost and prize collection profits in the objective function~\eqref{eq:obj}.
The value of $\rho$ can be interpreted as the maximum number of kilometres one is willing to drive to avoid an overflow, scaled by the estimated probability of overflow occurring: for example, with $\rho = 100$~\text{kilometres}, one would be willing to drive at most 100 kilometres to avoid a certain overflow, but only 50 kilometres if the estimated probability is 0.5.
The parameters of~$\epsilon$ and~$\rho$ require tuning for the specific setting of our case study.

The remainder of this section describes how to infer $\text{Pr}( Z_{c~i} > V_c)$, and how to estimate $Z_{c~i}$ from data.
Since inference and estimation are the same for each container cluster $c \in C$, we drop the dependency on $c$ for the remainder of this section.

\paragraph{Inference}
When planning a shift at time $T_{i}$, the total volume of deposits that arrives before the next planning time $T_{i + 1}$ is the given by the random variable
\begin{equation}
    \label{eq:Z_isr}
    Z_{i} = \sum_{j = N(A(T_i)) + 1}^{N(T_i)} v_{j} + \sum_{j = N(T_i) + 1}^{N(T_{i + 1})} v_{j}.
\end{equation} 
Clearly, the first term corresponds to the total volume of deposits known at time $T_i$.
The second term contains the volume from deposits arriving between $T_i$ and $T_{i+1}$.
Although the overflow probability of a cluster depends on the planned service time of a cluster rather than the next planning time $T_{i + 1}$, we neglect this dependency since the resulting service time-dependent prizes would significantly increase the complexity of the routing problem.
The numerical results of~\cref{sec:case_study} show that this simplification is not a grave error.
It follows from the central limit theorem that the first term is very well approximated by a normal random variable with distribution $\text{Norm}(n \mu, n \sigma^2)$.
For notational convenience, let $l_{i} = l(T_i, T_{i + 1})$ be the expected number of deposits arriving between $T_i$ and $T_{i + 1}$.
Then, as the second term is distributed as a compound Poisson random variable, it is well approximated by a normal random variable with distribution $\text{Norm}(l_i \mu, l_i \sigma^2 + l_i \mu^2)$.
As the sum of two normal random variables is again a normal random variable, we have that
\[ Z_{i} \sim \text{Norm}((n + l_i) \mu, (n + l_i) \sigma^2 + l_i \mu^2), \]
and thus
\begin{equation}
    \label{eq:overflow_estimate}
    \text{Pr}(Z_{i} > V) \approx 1 - \Phi \left( \frac{V - (n + l_i) \mu}{\sqrt{(n + l_i) \sigma^2 + l_i \mu^2)}} \right),
\end{equation}
where $\Phi(\cdot)$ is the cumulative distribution function of the standard normal distribution.
Thus, if we can estimate $\mu$ and $\sigma$ from data, the right-hand side of~\eqref{eq:overflow_estimate} allows us to compute the overflow probability based on the number of known deposits $n$ up to $T_i$ and the arrival rate $l_i$ until the next shift planning time $T_{i+1}$.

% We remark that the central limit theorem gives excellent approximations in our setting.
% To see this, note that the maximum volume per deposit is 100 litres, and the capacity of a container cluster is at least 5 m$^3$, so that (theoretically speaking) the container needs at least 50 deposits before it can overflow. 

\paragraph{Estimation}
It remains to estimate $\mu$ and $\sigma^2$.
Since the containers do not measure deposit volumes, we need an indirect method to estimate the mean $\mu$ and variance $\sigma^2$ of the deposit volume.
For this, we can count the number of deposits that take place between services, and rely on the drivers to register whether a container cluster has overflowed or not when they arrive at the cluster for service.
Thus, for each container cluster $c$, we observe a sequence of $(d, o)_j$ observations, where $d_j$ is the number of deposits that arrived between the $(j - 1)^\text{th}$ and $j^\text{th}$ service, and $o_j \in \{0, 1 \}$ indicates whether an overflow occurred.

Suppose we have a sample of $m$ observations at time $T_i$.
Write $P(\mu, \sigma \mid d) = \text{Pr}(S_d > V)$ for the overflow probability of the random variable $S_d \sim \text{Norm}(d\mu, d\sigma^2)$.
The likelihood of a single observation $(d_j, o_j)$ is
\[ L(\mu, \sigma \mid d_j, o_j) = P(\mu, \sigma \mid d_j)^{o_j} (1 - P(\mu, \sigma \mid d_j))^{1 - o_j}. \]
As the observations are assumed to be independent, the joint log likelihood of the entire sample can be written as
\[ \log L(\mu, \sigma \mid D, O) = \sum_{j = 1}^m \biggl[ o_j \log P(\mu, \sigma \mid d_j) + (1 - o_j) \log(1 - P(\mu, \sigma \mid d_j)) \biggr]. \]
We observed that solving $(\hat \mu, \hat \sigma) = \arg \max_{\mu, \sigma} \log L(\mu, \sigma \mid D, O)$  using standard gradient methods results in good estimates for $\mu$ and $\sigma^2$ despite $\log L(\mu, \sigma \mid D, O)$ not being concave.

Since we cannot guarantee finding an optimal $\hat \mu$ and $\hat \sigma$ given the sample data, we may instead want to choose $\hat \mu$ and $\hat \sigma$ in a conservative way, to ensure that we do not underestimate the overflow probabilities.
Now observe for a random variable $X$ with known mean $\mu$ and support $[0, b]$ that $X^{2} \leq b X$, hence it holds for the variance that $V(X) = E (X^2) - E(X)^2 \le E(bX) - \mu^2 = \mu(b - \mu)$.
Since the deposit volumes are necessarily bounded to the interval of $[0, 100]$ litres due to physical limitations of the container drums, this observation establishes a simple upper bound on $\sigma \le \sqrt{\mu(100 - \mu)}$.
Now, solving $\max_\mu \log L(\mu, \sqrt{\mu(100 - \mu)} \mid D, O)$ is a single-variable optimisation problem over the interval $(0, 100]$.
Although gradient methods need not find an optimal $\hat \mu$ since multiple locally optimal solutions might again exist, by starting the gradient method from the initial point $\hat \mu_0 = 100$ we can ensure the largest $\hat \mu$ is selected. 
The resulting values of $\hat \mu$ and $\hat \sigma = \sqrt{\hat \mu (100 - \hat \mu)}$ are conservative estimates for $\mu$ and $\sigma$.

\section{Groningen case study}
\label{sec:case_study}

In this section we solve the case study of the municipality of Groningen.
\Cref{subsec:simulator} explains the simulation environment we use to model the case study setting.
In~\cref{subsec:case_tuning} we tune our ISR policy's parameters.
We compare the ISR policy to the baseline policy in~\cref{subsec:results}.
In~\cref{subsec:sensors} and~\cref{subsec:fleet_size} we explore several additional scenarios that are of interest to the municipality.

\subsection{Simulation environment}
\label{subsec:simulator}

We implement the case study setting in a simulation environment which we describe in this section.
\makeatletter
\if@BLINDREV
The full software implementation of our simulation environment and policies is available at  \texttt{URL redacted for double blind review}.
\else
The full software implementation of our simulation environment and policies is available at \url{https://github.com/N-Wouda/Groningen-Waste-Collection}.
\fi
\makeatother

We obtained data from the municipality of Groningen about the container clusters (location, capacity, and volume correction factors), number of vehicles they currently operate, and all deposit arrivals at the container clusters during the first quarter of 2023.
The registered deposits amount to nearly two million data points, and allow us to estimate the rate parameters of the arrival process for each container cluster and hour-of-the-day.
We exclude five containers clusters from the data because they saw no deposits during the first quarter of 2023---these are all clusters that are either no longer, or not yet, in use.
% these include clusters that are not yet in use (2), used only for testing (1), a glass container not part of the residential waste collection system (1), and a temporary cluster no longer in use (1).
We sample the deposit volumes from a triangular distribution with minimum 10, mode 30, and maximum 60 litres.
The parameters of this distribution are loosely based on the results of a small-scale pilot study using fill level sensors that the municipality performed.
Given the locations of the depot and the container clusters, we estimate travel distances and durations between these locations from OpenStreetMap data through the OSRM project~\citep{luxen_real-time_2011}, using a heavy truck profile to model waste collection vehicles.

Our simulator is a discrete-event simulator, where different types of events are tracked in an event queue.
Events are handled in chronological order.
Over a fixed time horizon, the simulator inserts all deposit arrival and shift planning events in the queue before starting the simulation run.
These events are generated separately from the main simulation run, using a fixed seed.
This ensures our simulation runs benefit from common random numbers for variance reduction~\citep{glasserman_guidelines_1992}.

When a shift plan event is handled, the simulator uses one of the policies of~\cref{sec:policies} to generate service events and inserts these into the event queue.
Coffee and lunch break events are planned by the simulator (not the policy) as late as possible, such that the driver can still return to the depot for the break to start in time.
This ensures that policies that are not aware of the breaks (such as the baseline policy) do schedule them, although at sub-optimal positions in the route plan.
Once the simulation concludes, we process all events to compute the values of the performance measures described in~\cref{tab:performance_measures}.

\begin{table}
    \caption{Simulation performance measures.}
    \label{tab:performance_measures}

    \begin{threeparttable}
    \footnotesize

    \begin{tabular}{p{0.025\linewidth - 2\tabcolsep}p{0.3\linewidth-2\tabcolsep}p{0.675\linewidth-2\tabcolsep}}
        \toprule
        \multicolumn{2}{l}{Measure} & Definition \\
        \midrule
        \multicolumn{2}{l}{\textit{Operational costs}} & \\
        & Average daily driving distance & Average total distance (in kilometres). \\
        & Average route duration & Average route duration (in hours), including time spent on breaks. This performance measure indicates the amount of slack in each route, which is important to consider for future system growth. \\
        & Average number of routes & Average number of routes planned each day. Fewer routes equate to lower personnel costs and a smaller required number of collection vehicles. \\
        & Average number of clusters & Average number of clusters serviced each day. \\
        \midrule
        \multicolumn{2}{l}{\textit{Service quality}} & \\
        & Average service level & Percentage of services with a cluster fill level not exceeding 100\%. \\
        & Average cluster fill level & Average fill level percentage of each cluster upon service. \\
        & Average overflow volume & When a cluster has overflowed before service, this is the average amount of volume (in litres) that exceeds capacity at the service. This measure is useful to understand worst-case behaviour when a cluster is not serviced in time. \\
        & Number of unserviced clusters & Number of clusters that are never serviced over the simulation horizon. This measure helps understand system stability. \\
        \bottomrule
    \end{tabular}
    \end{threeparttable}
\end{table}

\subsection{Parameter tuning}
\label{subsec:case_tuning}

Our ISR policy requires us to tune $\epsilon$ in~\cref{eq:isr_prizes}, which controls when a container cluster must be serviced during the current shift, and $\rho$, which balances the distance cost and prize collection profits in the objective function~\eqref{eq:obj}.
Both parameters admit a helpful interpretation, which makes it relatively easy to determine a range of good potential values.
It makes sense to only require service when the model is certain the cluster will overflow before the next planning moment.
Thus, we should set $\epsilon$ small to ensure this is the case, and it seems reasonable from the context to consider $\epsilon \in (0, 0.01, \ldots, 0.1)$.
The parameter $\rho$ has the interpretation of a Lagrangian multiplier: $\rho$ is the maximum additional distance that may (profitably) be travelled to empty an overflowing cluster.
A reasonable estimate for $\rho$ would be to select either twice the average arc distance $\frac{2}{|A|} \sum_{(i, j) \in A} d_{ij}$, or twice the maximum distance $2 \max_{(i, j) \in A} d_{ij}$, since those values are of the same order of magnitude as the distances in the objective~\eqref{eq:obj}.
For Groningen, these numbers amount to 9 and 44 kilometres, respectively, suggesting that a range of values in $1, 2, 4, \ldots, 32, 64$ kilometres should contain a good parameter value for~$\rho$.
To further investigate the effects of larger $\rho$ on the service level, we additionally consider $\rho$ in $(128, 256, 512, 1024)$ kilometres.

We remark that we explicitly choose not to tune the routing solver's parameters for the case study setting.
Such tuning would undoubtedly result in better performance, especially with the limited runtime allowed for the routing solver inside the simulator.
However, this assumes the municipality has access to a routing solver they can tune, which, given commercial black box solvers like the one they currently use, is not guaranteed.
Second, it assumes they are willing to expend the (often significant) effort in terms of time and compute resources to tune a routing solver.
As such, we treat PyVRP as a black box routing solver in the context of our case study, and use its default parameter configuration.

We perform a grid search over all candidate pairs of $\epsilon$ and $\rho$.
To limit the computational burden we use a single run for each parameter pair.
This is acceptable since randomness does not significantly impact the outcomes obtained from our simulator, as we observed over multiple test runs.
Each daily routing problem is given only one minute of runtime to mimic real-time constraints.
We run the simulator for one year and three months for each of the parameter combinations, and discard the first three months of each simulation run as warm-up time.
We determined the warm-up time by visually inspecting plots of the performance measures over time, and observed that these stabilise already around one month into the simulation run.
We then added two more months as a safety margin before we began collecting results.
From the final year of simulated data we collect the average daily driving distance and the average service level as performance measures (see~\cref{tab:performance_measures} for definitions).

\begin{figure}
    \centering
    \includegraphics[width=\linewidth]{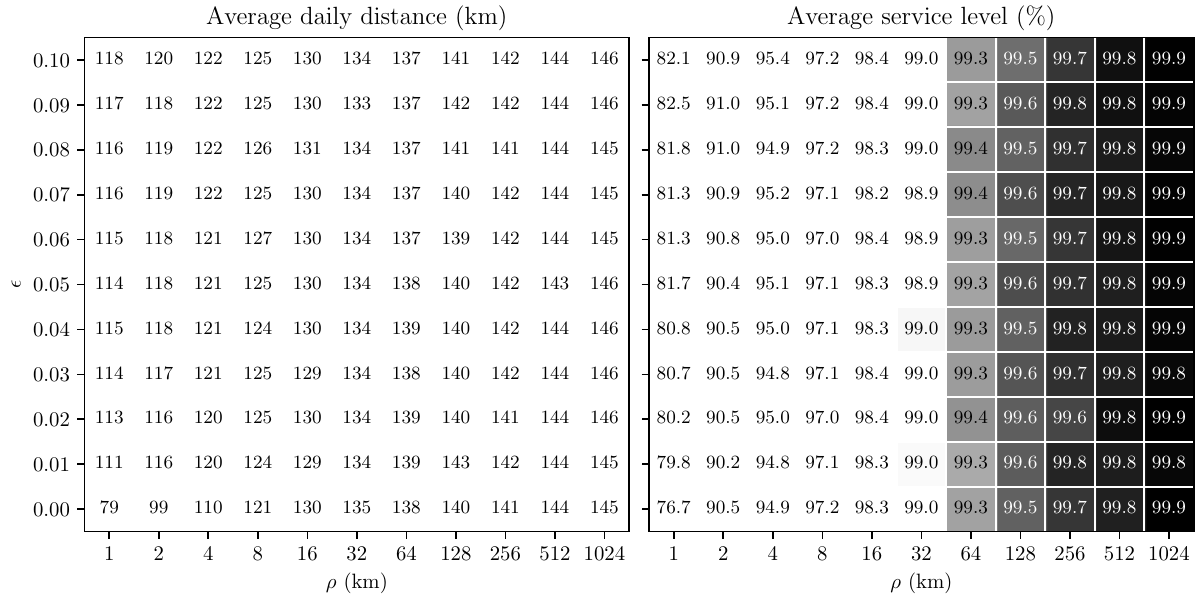}
    \caption{Parameter tuning results for different values of $\epsilon$ and $\rho$, for both the average daily driving distance (in kilometres) and the average service level (\%). Better results with respect to the average service level are shaded in darker colours.}
    \label{fig:case_tuning}
\end{figure}

\Cref{fig:case_tuning} shows the resulting average daily distance (in kilometres) and service levels.
When $\epsilon = 0$, that is, when there are no required services, we observe that the value of $\rho$ significantly impacts the results.
In particular, the average distance and service level, respectively, increase from 79 kilometres and 76.7\% when $\rho = 1$, to 145 kilometres and 99.9\% service level when $\rho = 1024$.
We conclude that $\rho$ is very effective at ensuring high service levels.
We observe that the effect of increasing $\rho$ tapers off after $\rho = 32$ kilometres, at which point a 99\% service level is first reached.
Remarkably, a near-100\% service level does not seem to require exorbitant operational costs: starting from 94.9\% (at $\rho = 4$), reaching 99.9\% (at $\rho = 1024$) increases the average daily driving distance by only 35 kilometres.

When $\epsilon > 0$, some cluster services may be required.
Such required services appear to be helpful in improving overall service levels when $\rho$ is small, but the effect quickly dissipates as $\rho$ increases.
For our case study, simply increasing $\rho$ thus appears to be much more effective than tuning $\epsilon$.
Based on~\cref{fig:case_tuning} we henceforth set $\epsilon = 0$ and just explore different values of $\rho$ to better understand the trade-off between routing costs and service level.

\subsection{Comparing the ISR and baseline policies}
\label{subsec:results}

We investigate the performance of the two policies of~\cref{sec:policies} to determine how well our ISR policy works compared to the baseline. 
% To this end, we again use the simulator, where we vary the policy used to solve the daily shift planning problem.
We derive the performance measures discussed in this section from the final year of simulated data, discarding the first three months of each simulation run as warm-up time.
Since the simulator draws container deposits randomly, each of the performance measures listed in~\cref{tab:performance_measures} and discussed in this section is averaged over ten simulation runs, each with different seeds.

We study different parametrisations of the two policies of~\cref{sec:policies} to determine the effect of different parameter choices.
In particular, for the baseline policy we investigate the effects of a different number of daily serviced clusters, ranging from 150 to 250 container clusters.
We will refer to the baseline policy with 250 clusters as the \textit{default} baseline policy, since this is the value the municipality currently uses.
For the ISR policy, we vary the scaling parameter $\rho$ from 4 kilometres to 1024 kilometres.
The results are given in~\cref{tab:simulation_results_4veh}.

\begin{table}
    \caption{%
        Simulation outcomes on the measures defined in~\cref{tab:performance_measures} for each policy, with different parametrisations, averaged over ten different seeds.
    }
    \label{tab:simulation_results_4veh}
    \centering

    \begin{threeparttable}
    \footnotesize

    \begin{tabular}{lrrrrrrrrr}
        \toprule
        & & \multicolumn{4}{l}{\textit{Operational costs}} & \multicolumn{4}{l}{\textit{Service quality}}  \\
        \cmidrule(lr){3-6}\cmidrule(lr){7-10}
        & & Distance & Duration & \# routes & \# clusters & Service level & Fill level & Overflow & \# unserviced  \\
        \midrule
        \multirow{6}{*}{\rotatebox[origin=c]{90}{\textit{Baseline}}} & \# clusters & & & & & & & & \\
        \cmidrule(lr){2-2}
        & 150 & 171.8 km & 5.6 h & 3.0 & 150 & 97.0\% & 68.1\% & 2492 L & 156 \\
        & 175 & 189.3 km & 6.5 h & 3.0 & 175 & 95.3\% & 69.9\% & 2845 L & 0 \\
        & 200 & 215.9 km & 7.5 h & 3.0 & 200 & 99.5\% & 61.2\% & 1077 L & 0 \\
        & 225 & 231.4 km & 6.1 h & 4.0 & 225 & 99.8\% & 54.4\% & 1210 L & 0 \\
        & 250 & 251.2 km & 6.8 h & 4.0 & 250 & 99.9\% & 48.9\% & 1537 L & 0 \\
        \midrule
        \multirow{6}{*}{\rotatebox[origin=c]{90}{\textit{ISR}}} & $\rho$ & & & & & & & & \\
        \cmidrule(lr){2-2}
        & 4 km    & 109.7 km & 3.2 h & 4.0 & 138.9 & 94.8\% & 87.9\% & 480 L & 1 \\
        & 16 km   & 129.8 km & 3.5 h & 4.0 & 143.6 & 98.3\% & 85.2\% & 130 L & 0 \\
        & 64 km   & 137.9 km & 3.6 h & 4.0 & 147.0 & 99.4\% & 83.2\% &  99 L & 0 \\
        & 256 km  & 142.8 km & 3.7 h & 4.0 & 150.1 & 99.7\% & 81.5\% &  90 L & 0 \\
        & 1024 km & 146.4 km & 3.8 h & 4.0 & 153.1 & 99.9\% & 79.9\% &  93 L & 0 \\
        \bottomrule
    \end{tabular}
    \end{threeparttable}
\end{table}

\paragraph{Baseline policy}
We observe that the default baseline policy servicing 250 clusters achieves an excellent performance level of 99.9\%, at a daily cost of on average 251.2 kilometres driving distance and an average shift duration of 6.8 hours.
To put this into perspective, a service level of 99.9\% indicates that an overflow due to late servicing of a container cluster is observed around once every four days, since $0.1\% \times 250 = 0.25$. 
However, although overflow events are quite rare under the baseline policy, when an overflow does occur, the average overflow volume is substantial: 1537 litres, or around one third of a container's worth of capacity.
Such a large overflow volume is not something the driver can quickly clean during their shift, and requires a clearner to be sent to the cluster.
Closer inspection of the simulation results reveals that the additional volume of deposits that causes overflow almost always arrives on the day of a planned service, with the cluster already being very full when the shift plan is made in the morning.

Under the default baseline policy, the average cluster fill level amounts to 48.9\%.
This observation corroborates an observation made by the municipality that clusters are not as full as they would like them to be, which they measure ad-hoc based on the container weight when it is lifted out of the ground during service.
Attempting to improve the cluster utilisation by reducing the number of serviced clusters each day helps, but not by much: reducing the number of clusters to 200 (from 250) increases the average fill level to 61.2\%, while the service level lowers to 99.5\%.
The average route duration increases to 7.5 hours, as the routing solver determines that routes with 200 clusters (barely) fit into the maximum shift duration of three vehicles.
In practice, however, a fourth vehicle is likely still needed to service 200 clusters a day.
The average fill level increases further as the number of serviced clusters drops below 200, but then the average service level drops significantly as well; at 175 clusters it is just 95.3\%.
Paradoxically, at 150 clusters the average service level is higher, at 97.0\%, but that is because now 156 clusters are never serviced at all because the service capacity of the vehicles is insufficient to keep up with all deposit arrivals.
From this, we conclude that using just the number of daily serviced containers does not offer much leeway to achieve better capacity utilisation, at least not without also sacrificing service level.

\paragraph{ISR policy}
Let us now turn to our ISR policy.
As we see in~\cref{fig:case_tuning}, the ISR policy with $\rho = 1024$ achieves the same service level as the baseline policy of 99.9\%, at significantly lower average daily driving distance: 146.4 kilometres, compared to 251.2 kilometres for the baseline policy.
The average daily driving distance is thus reduced by 105 kilometres, a reduction of 41\%.
Additionally, the ISR policy performs much better when an overflow event does occur: in this case, the average overflow is only 93 litres, as opposed to 1537 litres.
Since 93 litres equates to about two reasonably filled garbage bags, this is an amount of waste the driver might be able to clean up during the shift, which would significantly reduce the number of times a cleaner needs to be sent to the cluster.
Overall, the ISR policy better accounts for the arrival of new waste than the baseline policy, since average cluster fill levels are significantly higher (79.9\% vs. 48.9\%), while the same high service level is achieved and average overflow volume is much reduced.
Further, the average shift lasts only 3.8 hours as opposed to 6.8 hours under the default baseline policy, suggesting there is ample slack in the system to account for future expansion.
These observations lead to~\cref{in:isr_is_better}.
\begin{insight}
    \label{in:isr_is_better}
    Adopting our ISR policy parametrised with $\rho = 1024$ kilometres over the default baseline policy reduces driving distance by 41\% and route duration by 44\%, at no reduction in average service level.
\end{insight}

We further observe from~\cref{tab:simulation_results_4veh} that varying $\rho$ allows for fine-grained balancing of operational costs and service quality. 
Selecting a value of $\rho = 16~\text{kilometres}$ reduces the average distance by another 16 kilometres, shortens the average shift duration by another twenty minutes, while the service level is 98.3\%.
Such a service level may be too low in practice, but overflows are only 130 litres on average when $\rho$ is set to 16 kilometres.
Although the service level decreases, ten overflows of this form still amount to less waste on the streets than a single overflow under the default baseline policy.

\subsection{Effect of fill level sensors}
\label{subsec:sensors}

This section investigates the effect of fill level sensors that provide real-time information about the cluster fill levels.
This is motivated by a current pilot involving such sensors.
No decision has yet been reached on whether to invest in installing these sensors in all clusters.
For the investment to be worthwhile, substantial improvements in either operational cost or service level should be achieved.
As~\cref{subsec:results} shows that excellent service levels can already be reached without these sensors, we may only hope that operational costs can be decreased using the fill level information to make better servicing decisions.
We explore that here, by slightly modifying the policies to use fill level information when determining prizes, rather than the observed number of deposits.
The required policy modifications are presented in Appendix~\ref{app:sensors}.

\begin{table}
    \caption{%
        Simulation outcomes on the measures defined in~\cref{tab:performance_measures} for each policy, with and without sensors, averaged over ten different seeds.
    }
    \label{tab:simulation_results_sensors}
    \centering

    \begin{threeparttable}
    \footnotesize

    \begin{tabular}{lrrrrrrrr}
        \toprule
        Sensors? & & \multicolumn{3}{l}{\textit{Operational costs}} & \multicolumn{4}{l}{\textit{Service quality}}  \\
        \cmidrule(lr){3-5}\cmidrule(lr){6-9}
        & & Distance & Duration & \# clusters & Service level & Fill level & Overflow & \# unserviced  \\
        \midrule
        No & Baseline & 251.2 km & 6.8 h & 250 & 99.9\% & 48.9\% & 1537 L & 0 \\
        No & ISR & 146.4 km & 3.8 h & 153.1 & 99.9\% & 79.9\% &  93 L & 0 \\
        \midrule
        Yes & Baseline & 245.3 km & 6.8 h & 250 & 99.7\% & 48.9\% & 1614 L & 0 \\
        Yes & ISR & 143.6 km & 3.8 h & 150.4 & 99.9\% & 81.3\% &  81 L & 0 \\
        \bottomrule
    \end{tabular}
    \end{threeparttable}
\end{table}

\Cref{tab:simulation_results_sensors} presents the results of the default baseline policy and the ISR policy with $\rho = 1024$ kilometres, with and without fill level sensors.
We observe that fill-rate sensors marginally reduce the daily distance under both policies, by 5.9 and 2.8 kilometres for the baseline and ISR policy, respectively.
Perhaps surprisingly, the service level of the baseline policy worsens slightly, from 99.9\% to 99.7\%.
This is due to the estimate of $D$, the assumed deposit volume: the baseline policy using this assumption is more conservative than the fill level approach, and thus services clusters earlier.
On the whole, and perhaps surprisingly, the table suggests that while fill level sensors offer some operational benefits, these improvements are not substantial.
Given that such sensors are more expensive, can break down and require extra maintenance, it seems that equipping the containers with such sensors is not worth the cost.
This leads to~\cref{in:sensors_do_not_matter_much}.

\begin{insight}
    \label{in:sensors_do_not_matter_much}
    For both the baseline and ISR policies, the availability of fill level data from sensors offers only small improvements (less than 3\%) in driving distance, with no improvement in service level.
\end{insight}

% TODO sensors can also break down, are inaccurate, might not be worth

\subsection{Effect of fleet size}
\label{subsec:fleet_size}

Since there is significant slack in the shift plans using the ISR policy (as evidenced by the short route durations in~\cref{tab:simulation_results_4veh}), we investigate whether the fourth collection vehicle is actually required.
Reducing the number of collection vehicles is financially very attractive, as the purchasing cost lies between 500K\euro{} to 1M\euro{}, and operational and maintenance costs are quite significant.

\begin{table}
    \caption{%
        Simulation outcomes on the measures defined in~\cref{tab:performance_measures} for each policy using three vehicles, with different parametrisations, averaged over ten different seeds.
    }
    \label{tab:simulation_results_3veh}
    \centering

    \begin{threeparttable}
    \footnotesize

    \begin{tabular}{lrrrrrrrrr}
        \toprule
        & & \multicolumn{4}{l}{\textit{Operational costs}} & \multicolumn{4}{l}{\textit{Service quality}}  \\
        \cmidrule(lr){3-6}\cmidrule(lr){7-10}
        & & Distance & Duration & \# routes & \# clusters & Service level & Fill level & Overflow & \# unserviced  \\
        \midrule
        \multirow{4}{*}{\rotatebox[origin=c]{90}{\textit{Baseline}}} & \# clusters & & & & & & & & \\
        \cmidrule(lr){2-2}
        & 150 & 171.7 km & 5.6 h & 3.0 & 150 & 97.1\% & 68.1\% & 2464 L & 157 \\
        & 175 & 189.4 km & 6.5 h & 3.0 & 175 & 95.2\% & 69.9\% & 2797 L & 0 \\
        & 200 & 215.7 km & 7.5 h & 3.0 & 200 & 99.5\% & 61.2\% & 1101 L & 0 \\
        \midrule
        \multirow{6}{*}{\rotatebox[origin=c]{90}{\textit{ISR}}} & $\rho$ & & & & & & & & \\
        \cmidrule(lr){2-2}
        & 4 km    & 118.0 km & 4.6 h & 3.0 & 137.3 & 92.5\% & 88.8\% & 357 L & 3 \\
        & 16 km   & 139.2 km & 5.1 h & 3.0 & 142.5 & 97.3\% & 85.9\% & 143 L & 0 \\
        & 64 km   & 147.7 km & 5.3 h & 3.0 & 145.8 & 98.8\% & 83.9\% & 116 L & 0 \\
        & 256 km  & 152.2 km & 5.4 h & 3.0 & 149.0 & 99.4\% & 82.1\% & 107 L & 0 \\
        & 1024 km & 155.8 km & 5.5 h & 3.0 & 151.8 & 99.7\% & 80.6\% & 104 L & 0 \\
        \bottomrule
    \end{tabular}
    \end{threeparttable}
\end{table}

The results with three collection vehicles are given in~\cref{tab:simulation_results_3veh}.
We compare with the baseline policy with only three collection vehicles: in this case, servicing 200 or more clusters each day cannot feasibly be achieved within the shift time limit of 7 hours.
However, we observe that the baseline policy with fewer clusters does not work well: with 175 clusters the service level drops significantly (from 99.5\% to 95.2\%), while at 150 clusters the overall service capacity is insufficient to keep up with the arrival of deposits, as evidenced by the 157 unserviced clusters.

The ISR policy with three vehicles, however, appears to be almost as good as the four-vehicle policy.
In particular, at $\rho = 1024$, the average service level is 99.7\% (0.2 percentage points below the four-vehicle policy), for an average daily driving distance of 155.8 kilometres (up from 146.4 kilometres).
The increase in distance is largely due to the longer average route duration (5.5 hours rather than 3.8 hours), and the additional break that is then required because the shift no longer ends before the lunch break.
An average route duration of 5.5 hours leaves sufficient slack in the system, while all vehicles are utilised substantially better than they are under the four-vehicle policy.
These observations lead to~\cref{in:3_rather_than_4_vehicles}.
\begin{insight}
    \label{in:3_rather_than_4_vehicles}
    Our ISR policy with three vehicles performs almost as well as the four-vehicle policy (both with $\rho = 1024$ kilometres): its average distance is just 6\% higher, and its attained service level only 0.2 percentage points lower.
\end{insight}

We also briefly investigated a setting with just two vehicles, but found that this does not offer sufficient collection capacity to keep up with the rate of deposits.

% {
% \color{orange}
% TODO maybe some additional scenarios? Like:
% \begin{itemize}
%     \item
%     What if different clusters have different volume distributions? Right now we assume the same triangular distribution for each clusters. We know arrival rates are very different between clusters, why would that not also be the case for volumes?

%     \item
%     Also look at baseline w/ break planning?

%     \item 
%     Re-optimise during each break? This would not be super easy to do in the current simulator, but may be worthwhile.
% \end{itemize}
% }

\section{Practical considerations}
\label{sec:discussion}

\Cref{subsec:results} reveals that the ISR policy with a scaling parameter of $\rho = 1024$ kilometres outperforms the default baseline policy on all relevant (operational) performance measures.
It achieves the same high service level of 99.9\% while reducing the average daily distance by 41\% from 251.2 km to 146.4 km, and the average shift duration by 44\%, from 6.8 hours to 3.8 hours.
Additionally, the ISR policy significantly lowers the average overflow volume to approximately 100 litres, compared to more than 1500 litres under the default baseline policy.
Based on these observations, we recommend that the municipality adopt our ISR policy.

A consideration for the implementation of our policy is that the high service levels we observe in the simulator are unlikely to be reached in practice.
In particular, containers inside a cluster may jam due to technical failure, or a deposit getting stuck in the drum.
This results in waste ending up on the street that is not caused by overflow.
It is difficult to detect these jams using our data, and we have been unable to take technical failures into account in our simulation environment.
The municipality sometimes receives complaints about a jammed container, and then needs to send a mechanic to resolve the issue.
Alternatively, the municipality could opt to involve local residents, who may `adopt' their container cluster to receive keys to the containers and can help clean up small overflows or attempt to fix jammed containers.
This helps keep the local neighbourhood clean and tidy.
The process of `container adoption' has already been implemented successfully in, amongst others, the municipalities of Amsterdam and Rotterdam.

\paragraph{Strategic}
In addition to the operational performance of the ISR policy, in~\cref{subsec:sensors} we considered whether the municipality should install fill level sensors in each container, to obtain real-time information on fill levels.
Surprisingly,~\cref{in:sensors_do_not_matter_much} suggests that this additional information does not increase service levels, and decreases average daily distance by less than 3\%.
Also, sensors can break down, increasing the cost and complexity of container maintenance.
Since the average operational costs barely decrease, the fill level sensors do not seem to justify the expense.

An alternative modification that is more likely to be effective relates to cluster aggregation.
In this paper we look at whole clusters, not at the individual containers inside each cluster.
However, a deposit does not happen at the cluster level in practice: it happens at one of the containers that make up the cluster.
If a resident tries to open the drum of the first container, and it is full, they need to try the second, and so on until they find a container in the cluster that is not full to deposit their waste.
Rather than going through this process, a resident may instead decide to simply deposit their waste next to the containers.
We have not taken this behaviour into account in our simulator.
This behaviour occurs mostly in container clusters with more than two containers: in our case study, that describes only 25 out of 850 clusters.
The municipality may want to modify the containers in these 25 clusters to indicate visually whether they are full or not.

\paragraph{Tactical}
At the tactical level, the investigation of~\cref{subsec:fleet_size} suggests that our ISR policy is also effective with three collection vehicles, rather than the standard four.
This follows because the ISR policy better utilises the available cluster capacity than the baseline policy, as evidenced by the increased fill levels and reduced number of cluster visits per day.

However, as the containers have higher average fill levels before they are serviced, the collection vehicles fill up faster than they do currently.
This might cause vehicles to reach their capacities before the scheduled moments during the driver breaks when the vehicles are emptied.
If this turns out to be a problem, the municipality should take capacity constraints into account and include additional breaks or separate emptying moments.
The shorter route durations under the ISR policy should allow for an additional break, which could also have a secondary benefit of increasing work enjoyment for the drivers.

\section{Conclusion}
\label{sec:conclusion}

\paragraph{Summary}
This paper studies a real-world case of urban waste collection arising in the city of Groningen, in The Netherlands, where the municipality collects residential waste from underground containers, organised in container clusters.
Urban waste collection is a challenging topic in city operations.
We propose an efficient integrated selection and routing (ISR) policy that jointly optimises the cluster selection and routing decisions.
The integration is achieved by first estimating prizes that express the urgency of selecting a container cluster to empty, and then solving a prize-collecting vehicle routing problem to collect these prizes while minimising routing costs.
Through a large simulation study, we show that our ISR policy achieves excellent results: it maintains the high service levels the municipality finds important, while average driving distances reduce by more than 40\%.
Additionally, we show that it may be possible to reduce the number of collection vehicles from four to three vehicles, hence by 25\%.
Further, if overflows do occur, the ISR policy achieves much smaller overflow volumes.

\paragraph{Future work}
While we have shown in this paper that the municipality may be able to save substantial costs compared to their current baseline policy, we have not yet implemented the ISR policy in practice.
Piloting our ISR policy with the municipality of Groningen is thus a natural next step to take.
Additionally, many other urban municipalities in Western Europe organise residential waste collection in a similar manner.
They should be able to benefit from our ISR policy as well.

Methodologically, the value of our ISR policy lies in the integration of selection and routing decisions.
That approach appears to be generally useful in dynamic routing problems, including, for example, many forms of same-day delivery problems (see, for example, the recent works of~\cite{baty_combinatorial_2024} or~\cite{lan_iterative_2024}).
For this purpose, our prize-collecting extension to the open-source PyVRP routing solver should prove helpful.

% Acknowledgments here
\ACKNOWLEDGMENT{%
The authors thank Matthijs den Engelsman and Thomas Miedema of the municipality of Groningen for their expert insight and feedback.
Niels Wouda also thanks Leon Lan and Wouter Kool for helpful discussions and suggestions while implementing support for our case study setting in PyVRP.
}% Leave this (end of acknowledgment)

\clearpage
\begin{APPENDICES}
    \section{PyVRP extension details}
\label{app:pyvrp_details}

This appendix describes the changes we made to the metaheuristic in the PyVRP solver to enable it to solve the routing problems of this paper.
PyVRP~\citep{wouda_pyvrp_2024} internally uses a hybrid genetic search (HGS) algorithm that is derived from the HGS-CVRP algorithm of~\cite{vidal_hybrid_2022}.
We discuss relevant changes and some implementation details in~\cref{subsec:routing_solver}, and then investigate the performance of these changes on a new set of benchmark instances in~\cref{subsec:numerical_experiments}

\subsection{Solving the routing problem}
\label{subsec:routing_solver}

The HGS algorithm used inside PyVRP is an iterative genetic algorithm that maintains two separate subpopulations of feasible and infeasible solutions, respectively.
In each iteration, two parent solutions are selected from the population by means of a binary tournament.
These parent solutions are them combined using a crossover operator to generate an offspring solution.
The generated offspring solution is improved using a local search procedure, and the locally optimal offspring is then inserted into the appropriate subpopulation.
Once the subpopulation grows beyond its maximum size, a survivor selection procedure is triggered: this procedure removes duplicate solutions from the subpopulation, as well as any solutions with low fitness.
Fitness is computed as a function of the solution's objective value, and its diversity with respect to other solutions in the population.
Survivor selection thus ensures a balanced population of diverse and high quality solutions.
For our problem setting, we need to make several adjustments and extensions to the core HGS algorithm available in PyVRP, particularly to the local search.

During the local search phase of the algorithm, optional clients whose visits are not sufficiently attractive to remain in the solution are removed, and clients not currently in the solution whose visits result in a cost reduction are re-inserted.
Re-insertion of a client $i \in C$ happens by evaluating the best re-insert location in the granular neighbourhood $\mathcal{N}(i)$.
This granular neighbourhood is based on the concept of \textit{correlated clients}~\citep{vidal_hybrid_2013}, and explicitly takes prizes into account.
In particular, we compute the correlation measure $\gamma_{ij}$ for each arc $(i, j) \in A$ as
\[
    \gamma_{ij} = d_{ij} + \beta_\text{wait} \underbrace{(e_j - l_i - s_i - t_{ij})^+}_{\text{minimum wait time}} + \beta_\text{tw} \underbrace{(e_i + s_i + t_{ij} - l_j)^+}_{\text{minimum time warp}} - p_j,
\]
where $\beta_\text{wait}$ and $\beta_\text{tw}$ are parameters that control the impact of time-related aspects on the correlation measure, and $p_j$ is the prize that is obtained by visiting location $j \in C$ using arc $(i, j) \in A$.
Given these correlation values, we define a granular neighbourhood $\mathcal{N}(i)$ of size $k > 0$ for each client $i \in C$ by adding to $\mathcal{N}(i)$ the $k$ clients $j \in C \setminus \{ i \}$ whose arcs $(i, j)$ have the smallest correlation values $\gamma_{ij}$.
By construction, the granular neighbourhoods of each client do not include the depot, nor the clients themselves.

A re-insert move is evaluated based on its cost delta, which is the change in objective value of inserting the client in the best possible location in the granular neighbourhood.
Optional clients are re-inserted when the resulting cost delta is negative.
Note that required clients are always re-inserted.
Removal and re-insertion are handled before the main granular neighbourhood search on the current solution commences.
This ensures that the search operators which we discuss next only need to handle clients that are currently in the solution, which simplifies their implementation considerably.

To keep the local search conceptually simple, we apply only common search operators with proven value.
In particular, our local search consists of the following operators:
\begin{itemize}
    \item
    The $(n, m)$-\textsc{exchange} operators of~\cite{taillard_parallel_1993} for all combinations of $1 \le n \le 3$ and $0 \le m \le 3$ with $n \ge m$, which swap $n$ clients from one route with $m$ clients from the other.
    This operator includes the well-known \textsc{relocate} ($n=1$, $m=0$) and \textsc{swap} ($n = m = 1$) operators as special cases, but generalises those to consider larger sequences of clients.

    \item
    The \textsc{tail} operator of~\cite{accorsi_fast_2021}, which swaps the tails of two routes by breaking the two routes in half just after two given clients.
    This operator is also known as \textsc{2-opt*} in the vehicle routing literature~\citep{toth_vehicle_2002}.

    \item
    The \textsc{swap*} operator of~\cite{vidal_hybrid_2022}, which considers the best swap move between two routes, but does not require that the swapped clients are inserted in each others place.
    Instead, each is inserted into the best location in the other route.
\end{itemize}

\subsection{Numerical experiments}
\label{subsec:numerical_experiments}

This section analyses the performance of our prize-collecting metaheuristic (as implemented in PyVRP version 0.7.0) on a set of large, new benchmark instances for the PCVRPTW.
We focus on the PCVRPTW specifically to enable a comparison with the literature: this problem variant captures the core difficulty optional clients represent, without requiring the additional details our real-world problem setting also needs, but which complicate comparisons with the existing literature.
We compare our implementation against best-known solutions from many earlier and longer runs, and against a strong baseline algorithm from the literature.
Each instance is solved in a single hour of runtime, using a single core of an AMD EPYC 7763 2.45GHz CPU.

\subsubsection{Benchmark instances}
\label{subsubsec:benchmark_instances}

We generate several instances to benchmark our prize-collecting metaheuristic.
Since we are primarily interested in performance on large instances of size similar to the case study, we take the well-known VRPTW instances of~\cite{gehring_parallel_1999} with 1000 clients as a basis.
These sixty instances have clients with demand $q_i > 0$, for $i \in C$, and a homogeneous vehicle capacity of $Q$. 
We add prizes on top of these instances, in the manner of~\cite{bulhoes_vehicle_2018}: for each client $i \in C$, we generate a prize $p_i = \max \{ h_i q_i, 1 \}$, where $h_i$ is sampled i.i.d. from $U[0.75, 2.25]$.
In expectation this results in prizes one-and-a-half times client demand, which is in line with the prize-generating procedure of~\cite{stenger_prize-collecting_2013}.
These generated instances and the latest best-known solutions are available at \url{https://github.com/PyVRP/Instances}.

\subsubsection{Computational results}

We compare our implementation in PyVRP against the recent prize-collecting HGS (PC-HGS) algorithm of~\cite{baty_combinatorial_2024}.
Since we have access to PC-HGS's implementation, we run both algorithms on the same hardware with the same stopping criterion of a single hour of run-time.
We use both algorithms with their default parameter configurations.

\begin{figure}
    \centering
    \includegraphics[width=\linewidth]{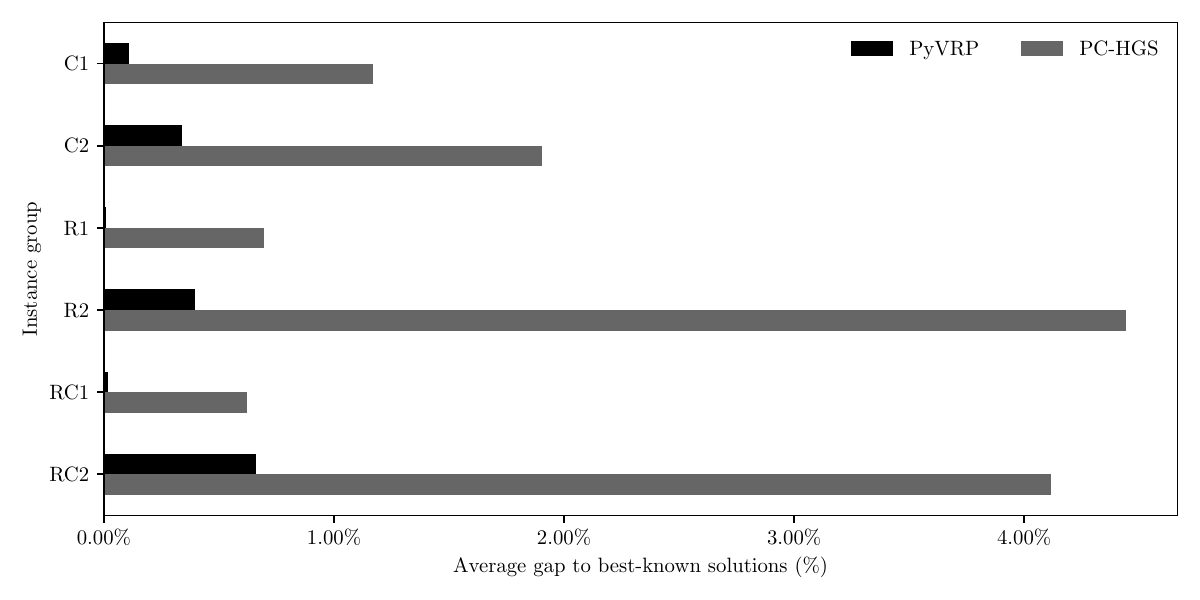}
    \caption{Average gaps to best-known solutions on different PCVRPTW instance groups after one hour of run-time, over ten different seeds, for both our implementation in PyVRP and the PC-HGS algorithm of~\cite{baty_combinatorial_2024}.}
    \label{fig:pyvrp_kleopatra_gaps}
\end{figure}

The average gaps to the best-known solutions (as of 24 March 2024) over ten different seeds are presented in~\cref{fig:pyvrp_kleopatra_gaps}, with details per instance presented in~\cref{tab:app_benchmark}.
We break the gaps down by instance group in the manner of~\cite{solomon_algorithms_1987}: groups C are clustered, R random, and RC semi-clustered; groups of type 1 have tight time windows and small vehicle capacities, whereas groups of type 2 have larger time windows and vehicle capacities.
Our implementation achieves an average gap of only 0.26\%, whereas PC-HGS has an average gap of 2.16\%.
Instances of type 2 appear to be particularly difficult, as both algorithms achieve their largest gaps on these instance groups.
This is likely due to the increased freedom the large time windows and vehicle capacities offer, which makes the search space for good, feasible solutions much larger.

{
\LTcapwidth=0.8\textwidth
\setlength{\tabcolsep}{8pt}
\begin{longtable}{lrrrrrrr}
\caption{
    Full benchmark results on the PCVRPTW instances over ten different seeds and one hour of runtime, for our implementation in PyVRP and the PC-HGS algorithm of~\cite{baty_combinatorial_2024}.
    We present the costs of the best-known solution (as of 24 March 2024), the average and best costs obtained by both algorithms, and the average percentage gaps.
    The best average gap is highlighted for each instance.
}
\label{tab:app_benchmark} \\

\toprule
Instance & BKS & \multicolumn{3}{c}{PyVRP} & \multicolumn{3}{c}{PC-HGS} \\
         &     & Avg.  & Best & Avg. gap   & Avg.   & Best & Avg. gap   \\
\cmidrule(lr){1-1}\cmidrule(lr){2-2}\cmidrule(lr){3-5}\cmidrule(lr){6-8}
\endfirsthead
\toprule
Instance & BKS & \multicolumn{3}{c}{PyVRP} & \multicolumn{3}{c}{PC-HGS} \\
         &     & Avg. & Best  & Avg. gap   & Avg. & Best   & Avg. gap   \\
\cmidrule(lr){1-1}\cmidrule(lr){2-2}\cmidrule(lr){3-5}\cmidrule(lr){6-8}
\endhead
\midrule
\multicolumn{8}{r}{Continued on next page} \\
\midrule
\endfoot
\bottomrule
\endlastfoot
\texttt{C1\_10\_1} & 24542.1 & 24545.7 & 24539.1 & \textbf{0.01\%} & 24948.0 & 24683.6 & 1.65\% \\
\texttt{C1\_10\_2} & 24901.5 & 24920.8 & 24901.5 & \textbf{0.08\%} & 25158.7 & 25104.8 & 1.03\% \\
\texttt{C1\_10\_3} & 24502.8 & 24560.3 & 24502.8 & \textbf{0.23\%} & 24720.1 & 24653.3 & 0.89\% \\
\texttt{C1\_10\_4} & 24456.2 & 24590.4 & 24456.2 & \textbf{0.55\%} & 24614.6 & 24542.4 & 0.65\% \\
\texttt{C1\_10\_5} & 24862.5 & 24880.7 & 24860.0 & \textbf{0.07\%} & 25297.5 & 25099.8 & 1.75\% \\
\texttt{C1\_10\_6} & 24676.3 & 24699.9 & 24676.3 & \textbf{0.10\%} & 25069.9 & 24851.0 & 1.59\% \\
\texttt{C1\_10\_7} & 24470.1 & 24464.7 & 24463.1 & \textbf{-0.02\%} & 24815.2 & 24719.9 & 1.41\% \\
\texttt{C1\_10\_8} & 24430.5 & 24433.2 & 24427.3 & \textbf{0.01\%} & 24734.5 & 24663.5 & 1.24\% \\
\texttt{C1\_10\_9} & 24589.1 & 24601.3 & 24589.1 & \textbf{0.05\%} & 24819.0 & 24734.1 & 0.94\% \\
\texttt{C1\_10\_10} & 24864.9 & 24867.9 & 24860.5 & \textbf{0.01\%} & 25003.8 & 24927.5 & 0.56\% \\
\texttt{C2\_10\_1} & 16585.5 & 16611.3 & 16589.2 & 0.16\% & 16595.7 & 16581.6 & \textbf{0.06\%} \\
\texttt{C2\_10\_2} & 16212.8 & 16250.3 & 16212.8 & \textbf{0.23\%} & 16392.0 & 16272.0 & 1.11\% \\
\texttt{C2\_10\_3} & 15718.4 & 15762.9 & 15723.7 & \textbf{0.28\%} & 16468.3 & 16168.1 & 4.77\% \\
\texttt{C2\_10\_4} & 15096.9 & 15201.0 & 15096.9 & \textbf{0.69\%} & 15803.0 & 15681.0 & 4.68\% \\
\texttt{C2\_10\_5} & 16257.1 & 16273.5 & 16261.0 & \textbf{0.10\%} & 16335.6 & 16285.6 & 0.48\% \\
\texttt{C2\_10\_6} & 16004.3 & 16061.4 & 16008.9 & \textbf{0.36\%} & 16228.6 & 16071.9 & 1.40\% \\
\texttt{C2\_10\_7} & 16017.7 & 16064.5 & 16020.5 & \textbf{0.29\%} & 16257.5 & 16073.7 & 1.50\% \\
\texttt{C2\_10\_8} & 15706.7 & 15769.0 & 15730.0 & \textbf{0.40\%} & 15949.3 & 15780.8 & 1.54\% \\
\texttt{C2\_10\_9} & 15819.2 & 15894.5 & 15830.8 & \textbf{0.48\%} & 16028.1 & 15952.1 & 1.32\% \\
\texttt{C2\_10\_10} & 15362.5 & 15425.5 & 15365.2 & \textbf{0.41\%} & 15698.5 & 15589.4 & 2.19\% \\
\texttt{R1\_10\_1} & 26273.5 & 26273.4 & 26271.1 & \textbf{-0.00\%} & 26498.8 & 26444.6 & 0.86\% \\
\texttt{R1\_10\_2} & 26312.5 & 26315.3 & 26312.5 & \textbf{0.01\%} & 26487.2 & 26453.2 & 0.66\% \\
\texttt{R1\_10\_3} & 25615.7 & 25618.2 & 25615.7 & \textbf{0.01\%} & 25789.6 & 25755.1 & 0.68\% \\
\texttt{R1\_10\_4} & 25212.2 & 25217.8 & 25212.2 & \textbf{0.02\%} & 25373.5 & 25333.3 & 0.64\% \\
\texttt{R1\_10\_5} & 25788.2 & 25788.4 & 25788.2 & \textbf{0.00\%} & 25991.7 & 25959.5 & 0.79\% \\
\texttt{R1\_10\_6} & 25397.0 & 25399.3 & 25397.0 & \textbf{0.01\%} & 25548.8 & 25517.8 & 0.60\% \\
\texttt{R1\_10\_7} & 25137.8 & 25138.2 & 25137.8 & \textbf{0.00\%} & 25317.3 & 25294.9 & 0.71\% \\
\texttt{R1\_10\_8} & 24896.0 & 24902.2 & 24896.0 & \textbf{0.02\%} & 25019.5 & 24996.2 & 0.50\% \\
\texttt{R1\_10\_9} & 25433.6 & 25434.1 & 25433.4 & \textbf{0.00\%} & 25618.5 & 25558.3 & 0.73\% \\
\texttt{R1\_10\_10} & 25791.5 & 25790.8 & 25790.3 & \textbf{-0.00\%} & 25995.2 & 25960.1 & 0.79\% \\
\texttt{R2\_10\_1} & 24173.0 & 24180.8 & 24115.9 & \textbf{0.03\%} & 24691.8 & 24472.4 & 2.15\% \\
\texttt{R2\_10\_2} & 20975.8 & 21067.2 & 20998.0 & \textbf{0.44\%} & 21992.5 & 21686.8 & 4.85\% \\
\texttt{R2\_10\_3} & 18269.5 & 18378.5 & 18278.3 & \textbf{0.60\%} & 19372.0 & 18992.5 & 6.03\% \\
\texttt{R2\_10\_4} & 15689.5 & 15764.2 & 15699.1 & \textbf{0.48\%} & 16569.0 & 16422.5 & 5.61\% \\
\texttt{R2\_10\_5} & 23465.2 & 23512.6 & 23461.4 & \textbf{0.20\%} & 24133.9 & 23791.2 & 2.85\% \\
\texttt{R2\_10\_6} & 20423.1 & 20512.4 & 20457.2 & \textbf{0.44\%} & 21518.1 & 21260.5 & 5.36\% \\
\texttt{R2\_10\_7} & 17790.1 & 17938.7 & 17795.3 & \textbf{0.84\%} & 18772.6 & 18437.6 & 5.52\% \\
\texttt{R2\_10\_8} & 15542.4 & 15619.2 & 15543.3 & \textbf{0.49\%} & 16305.2 & 16014.2 & 4.91\% \\
\texttt{R2\_10\_9} & 23016.3 & 23060.4 & 22990.7 & \textbf{0.19\%} & 23705.1 & 23414.9 & 2.99\% \\
\texttt{R2\_10\_10} & 22133.0 & 22189.9 & 22133.0 & \textbf{0.26\%} & 23056.9 & 22746.4 & 4.17\% \\
\texttt{RC1\_10\_1} & 24816.2 & 24817.2 & 24816.1 & \textbf{0.00\%} & 24926.3 & 24881.6 & 0.44\% \\
\texttt{RC1\_10\_2} & 25043.1 & 25048.5 & 25043.1 & \textbf{0.02\%} & 25226.3 & 25154.2 & 0.73\% \\
\texttt{RC1\_10\_3} & 24461.4 & 24468.4 & 24461.4 & \textbf{0.03\%} & 24652.4 & 24595.3 & 0.78\% \\
\texttt{RC1\_10\_4} & 24495.6 & 24503.4 & 24495.6 & \textbf{0.03\%} & 24629.1 & 24598.0 & 0.55\% \\
\texttt{RC1\_10\_5} & 25114.2 & 25115.0 & 25102.6 & \textbf{0.00\%} & 25214.3 & 25178.6 & 0.40\% \\
\texttt{RC1\_10\_6} & 24627.8 & 24633.5 & 24623.5 & \textbf{0.02\%} & 24832.4 & 24797.2 & 0.83\% \\
\texttt{RC1\_10\_7} & 24874.2 & 24881.3 & 24874.2 & \textbf{0.03\%} & 24987.9 & 24955.9 & 0.46\% \\
\texttt{RC1\_10\_8} & 24549.7 & 24556.3 & 24549.7 & \textbf{0.03\%} & 24642.7 & 24622.1 & 0.38\% \\
\texttt{RC1\_10\_9} & 24435.3 & 24437.6 & 24435.3 & \textbf{0.01\%} & 24679.6 & 24605.7 & 1.00\% \\
\texttt{RC1\_10\_10} & 24616.4 & 24618.6 & 24615.2 & \textbf{0.01\%} & 24781.6 & 24707.6 & 0.67\% \\
\texttt{RC2\_10\_1} & 19726.8 & 19830.8 & 19752.9 & \textbf{0.53\%} & 20545.4 & 20266.0 & 4.15\% \\
\texttt{RC2\_10\_2} & 17359.1 & 17590.1 & 17364.5 & \textbf{1.33\%} & 18565.9 & 18341.3 & 6.95\% \\
\texttt{RC2\_10\_3} & 15419.6 & 15564.3 & 15432.9 & \textbf{0.94\%} & 16229.0 & 16061.3 & 5.25\% \\
\texttt{RC2\_10\_4} & 13917.6 & 14043.4 & 13921.4 & \textbf{0.90\%} & 14652.4 & 14407.7 & 5.28\% \\
\texttt{RC2\_10\_5} & 18198.7 & 18303.5 & 18198.7 & \textbf{0.58\%} & 18831.0 & 18703.4 & 3.47\% \\
\texttt{RC2\_10\_6} & 18257.0 & 18319.0 & 18257.0 & \textbf{0.34\%} & 18949.2 & 18765.7 & 3.79\% \\
\texttt{RC2\_10\_7} & 17718.4 & 17786.5 & 17718.4 & \textbf{0.38\%} & 18349.1 & 18039.6 & 3.56\% \\
\texttt{RC2\_10\_8} & 17046.5 & 17114.5 & 17075.8 & \textbf{0.40\%} & 17567.0 & 17321.3 & 3.05\% \\
\texttt{RC2\_10\_9} & 16576.0 & 16680.3 & 16636.2 & \textbf{0.63\%} & 17033.5 & 16917.7 & 2.76\% \\
\texttt{RC2\_10\_10} & 16625.9 & 16721.3 & 16625.9 & \textbf{0.57\%} & 17109.1 & 16930.4 & 2.91\% \\
\end{longtable}
}

    \section{Policy modifications for fill level sensors}
\label{app:sensors}

We will present only the required modifications to the policies of~\cref{sec:policies}, in order of presentation, beginning with the changes to the baseline policy before turning to the ISR policy.
In addition to the cluster capacity $V_c$, assume we are now also given a used capacity $U_c$.

\subsection{Baseline policy}
We first recall~\eqref{eq:time_till_full_baseline}:
\[
    h_c = \max \left\{ t~\middle|~l_c(T_i, t) \le r_c \left( V_c / D - n_c \right) \right\},
\]
Here, the part requiring modification is given by $V_c / D - n_c$, which is the assumed number of additional deposits cluster $c \in C$ can still receive before overflowing.
We modify this to account for the actual remaining unused capacity $V_c - U_c$, divided by the assumed deposit volume $D$, to arrive at
\[
    h_c = \max \left\{ t~\middle|~l_c(T_i, t) \le r_c \left( \frac{V_c - U_c}{D} \right) \right\}.
\]
With this value for $h_c$, the rest of the policy proceeds as before.

\subsection{ISR policy}
We again suppress the dependency on $c$ to ease the notation.
Recall the definition of $Z_i$ in~\eqref{eq:Z_isr}:
\[
    Z_{i} = \sum_{j = N(A(T_i)) + 1}^{N(T_i)} v_{j} + \sum_{j = N(T_i) + 1}^{N(T_{i + 1})} v_{j},
\]
where the first term lists the volume of deposits already in the cluster, and the second term lists deposits that will occur between the current shift planning moment $T_i$ and the next planning moment at $T_{i + 1}$.
The first term is simply $U$, and no longer uncertain.
Thus, after applying the central limit theorem to the second term, we have that now
\[ Z_{i} \sim \text{Norm}(U + l_i \mu, l_i \sigma^2 + l_i \mu^2), \]
and thus
\[
    \text{Pr}(Z_{i} > V) \approx 1 - \Phi \left( \frac{V - U - l_i \mu}{\sqrt{l_i \sigma^2 + l_i \mu^2)}} \right).
\]
The rest of the policy proceeds as before.

\end{APPENDICES}

\clearpage
\bibliographystyle{informs2014trsc}
\bibliography{references} 

\end{document}